%AMSTeX File

% \masl13.tex --- 29-12-95 13:04

\documentstyle{amsppt}
\TagsOnRight
\NoBlackBoxes
\hoffset=.1in

\define\diag{\operatorname{diag}}
\define\sh{\operatorname{sh}}

\define\ld{\ldots}

\define\a{\alpha}
\define\be{\beta}

\define\de{\delta}

\define\vd{\varDelta}

\define\pti{\widetilde P}
\define\mti{\widetilde M}

\define\pd#1#2{\dfrac{\partial#1}{\partial#2}}

\define\vc#1{(#1_1,\ldots,#1_n)}
\define\vct#1{[#1_1,\ldots,#1_n]}
\define\vect#1{\{#1_1,\ldots,#1_n\}}

%\define\sumka{\sum^{}_{}}    %szerokie wskazniki

\define\ft{\infty}

\define\jak{\dfrac{i}{\kappa}}
\define\jakx{\dfrac{1}{\kappa}}
\define\jakkx{\dfrac{1}{\kappa^2}}

\define\leg{e^{-g\slash \kappa}}
\define\legg{e^{-2g\slash \kappa}}
\define\lefik{e^{-\pti_0\slash \kappa}}
\define\lef{e^{\pti_0\slash \kappa}}
\define\lefx{e^{\pti_0\slash 2\kappa}}
\define\leff{e^{-2\pti_0\slash \kappa}}
\define\ika{\dfrac{i\kappa}{2}}
\define\kdwa{\dfrac{\kappa}{2}}
\define\dekx{\dfrac{i}{2\kappa}} 
\define\kax{\dfrac{1}{2\kappa}}
\define\kapti{\Big(\dfrac{\widetilde{P}_0}{\kappa}\Big)}
\define\kaptid{\Big(\dfrac{\widetilde{P}_0}{2\kappa}\Big)}

\define\1{Poincar\'e algebra}
\define\2{momentum}
\define\3{dimensional}
\define\4{algebra}
\define\5{nonlinear}
\define\6{representation}
\define\7{commut}
\define\8{$\kappa$-Poincar\'e algebra}
\define\9{$\kappa$-Weyl algebra}
\define\0{deformation}

\document

\topmatter
\title Deformation map for generalized $\kappa$-Poincar\'e and $\kappa$-Weyl 
algebras
\endtitle
\rightheadtext{Deformation map} %{$\kappa$-Poincar\'e group}
\author Stefan Giller$^*$, Cezary Gonera$^*$, Micha\l \/ Majewski$^*$\\
{\it Department of Theoretical Physics}\\
{\it University of \L \'od\'z}\\
{\it ul. Pomorska 149/153, 90--236 \L \'od\'z, Poland}
\endauthor
\leftheadtext{S. Giller, C. Gonera, M. Majewski}
\thanks
*\ Supported by KBN grant 2\,P\,3022\,1706\,p\,02
\endthanks

\abstract A nonlinear transformation in the momentum space is constructed 
which converts the deformed action of Lorentz and Weyl generators on momenta 
into the standard one.
\endabstract

\endtopmatter
\document
\head 1. Introduction
\endhead

The $\kappa$-Poincar\'e algebra was introduced by Lukierski, Nowicki and 
Ruegg [1], [2]. This Hopf algebra deformation of classical Poincar\'e 
algebra seems to be quite interesting since it depends on dimensionful 
parameter $\kappa$. Many formal properties of \8 were analysed in a  number 
of recent paper. In particular, Majid and Ruegg [3] have shown that it 
possesses a biocrossproduct structure, one of the factors being the 
classical Lorentz \4 while remaining one --- a deformed (in co\4 sector) \2 
\4. In the \4ic sector the  whole \0 consists in \5 action of  Lorentz 
generators on momenta. It appears ([2], [4]) that one can find a \5 change 
of variables in \2 space which reduces the 'deformed' action of Lorentz 
generators to the standard one. The existence of such a mapping ('\0 map') 
does not mean that the deformed \4 is equivalnet to the standard one --- the 
co\4ic sectors are not mapped onto each other.

Quite recently([5]), the $\kappa$-deformed \1 has been generalized to the 
case of inhomogeneous \4 acting in arbitrary $n$-\3 flat space. We do not 
assume any longer that $g_{\mu\nu} = \diag (+,-,\ld,-)$; on the contrary, 
$g_{\mu\nu}$ is now arbitrary symmetry invertible matrix.

The resulting \4 reads
$$
\aligned
& [M^{\mu\nu}, M^{\a\be}]  = i(g^{\mu\be } M^{\nu\a} - g^{\nu\be} M^{\mu\a} 
+ g^{\nu\a} M^{\mu\be} - g^{\mu\a} M^{\nu\be}),\\
& [\pti_{\mu},\pti_{\nu}] = 0,\\
& [M^{ij},\pti_{0}] = 0,\\
& [M^{ij},\pti_{k}] = i\kappa(\de^j{}_kg^{0i} - \de^i{}_kg^{0j})(1 - \lefik) 
+ i(\de^j{}_k g^{is} - \de^i{}_k g^{js}) \pti_s,\\
& [M^{i0},\pti_{0}]= i\kappa g^{i0} (1 - \lefik) + i  g^{ik} \pti_k,\\
& [M^{i0},\pti_{k}] = - \dfrac{i\kappa}{2} g^{00}  \de^i{}_k (1 - \leff)  - 
i \de^i{}_k g^{0s} \pti_s \lefik \\
&  \phantom{[M^{i0},\pti_{k}] =} + i g^{0i} \pti_k ( \lefik -1) + 
\dfrac{i}{2\kappa} \de^i{}_k  g^{rs} \pti_r \pti_s - \jak g^{is} \pti_s 
\pti_k 
\endaligned
\tag"{(1a)}"
$$
and
$$
\aligned
& \vd \pti_0 = I \otimes \pti_0 + \pti_0 \otimes I,\\
& \vd \pti_k =  \pti_k \otimes \lefik + I \otimes \pti_k,\\
& \vd M^{ij} = M^{ij}  \otimes I + I \otimes M^{ij},\\
& \vd M^{i0} =  I \otimes M^{i0} + M^{i0} \otimes \lefik - \jakx M^{ij} 
\otimes \pti_j
\endaligned
\tag"{(1b)}"
$$
with $i,j,k,\ld = 1,2,\ld,n-1$.

It was also shown ([5]) that, provided $g_{00} = 0$ (which excludes the 
positive-definite metric $g_{\mu\nu}$), \4 (1) can be extended to 
$\kappa$-deformed Weyl \4. The additional operator $D$ (dilatation)  \7us 
with $M_{\mu\nu}$ and obeys
$$
\aligned
& [D,\pti_0] = i \kappa(1 - \lefik),\\
& [D,\pti_i] = i \pti_i  \lefik + i g_{0i} g^{0s} \pti_s (1 - \lefik) + 
\dekx g_{0i} g^{rs} \pti_r \pti_s\\
&  \phantom{[D,\pti_i] =} + \ika  g^{00} g_{i0} ( 1 -\lefik)^2,\\
& \vd D = D \otimes I + I \otimes D + g_{0i} M^{i0} \otimes ( 1 -\lefik) - 
\jakx g_{0i} M^{ik} \otimes \pti_k.
\endaligned
\tag{2}
$$

In the present paper we extend the construction of \0 map given in [4] to 
cover  \4s (1) and (2). Such a \0 map is very useful in classifying the \6s 
of deformed \4s.  Apart from this it provides an independent consistency 
check (of \4ic sector of) of \4s (1) and (2).

\head II. Deformation map
\endhead

In this section we will generalize to the case of general \8 as well as \9 
the \0 map obtained in [4]. To this end let $P_\mu$ and $M^{\mu\nu}$ be the 
generators of classical \1. We will be looking for the  \0 map of the 
following form
$$
\aligned
& \pti_0 = g(P_0,M^2),\\
& \pti_k = f(P_0,M^2) P_k + g_{k0} h(P_0,M^2);
\endaligned
\tag{3}
$$
here $M^2$ is the standard mass squared Casimir, $M^2 = g^{\mu\nu} 
P_{\mu}P_{\nu}$. Inserting (3) into (1a) we  arrive at the following set of 
equations (primes denotes differentiation with respect to $P_0$)
$$
\aligned
& g' P_0 = \kappa (1 - \leg) - g_{00}h,\\
& g' = f,\\
& f P_0 = \kappa (1 - \leg) - g_{00}h,\\
& f' P_0 = \jakx fh g_{00} + (\leg - 1) f,\\
& f' = - \jakx f^2,\\
&  \jakx fh = -h',\\
& f P_0 = \kdwa(1 - \legg) - g_{00}h \leg - \kax g_{00}h^2 + \kax P^2_0 
f^2,\\
& f = \leg f + \jakx  g_{00} fh + \jakx  P_0 f^2,\\
& h' P_0 = h (\leg - 1) + \jakx g_{00} h^2,\\
& h\leg  = \kax M^2 f^2 - \kax g_{00} h^2.
\endaligned
\tag{4}
$$

These equations can be solved to yield
$$
\aligned
& f = \dfrac{\kappa}{P_0 + C(M^2)} ,\\
& h = \dfrac{\kappa A(M^2)}{P_0 + C(M^2)} ,\\
& g = \kappa \ln \Big( \dfrac{P_0 + C(M^2)}{C - g_{00} A(M^2)} \Big)
\endaligned
\tag{5}
$$
with $A$ and $C$ subjected to the condition
$$
 g_{00} A^2(M^2) - 2A(M^2)  C(M^2) + M^2 = 0.
\tag{6}
$$
Finally, the \0 map reads
$$
\aligned
& \pti_0 =  \kappa \ln \Big( \dfrac{P_0 + C}{C - g_{00} A} \Big),\\
& \pti_i  = \dfrac{\kappa P_i}{P_0 + C}  + \dfrac{\kappa A}{P_0 + C} g_{i0}.
\endaligned
\tag{7}
$$
The inverse map takes the form
$$
\aligned
& P_0 =  (C - g_{00} A) \lef - C,\\
& P_i  = \dfrac{C - g_{00} A}{\kappa} \lef \pti_i -  g_{i0} A.
\endaligned
\tag{8}
$$
Let us insert these expressions into the formula for the Casimir operator 
$M^2$. After some \4 we obtain
$$
\aligned
\dfrac{2A}{C - g_{00} A} & =  \jakkx \Big[g^{00}\Big(2\kappa \sh 
\kapti\Big)^2 + 4\kappa g^{0l} \pti_l \lefx \sh \kaptid\\
& + g^{rs} \pti_r \lefx \pti_s \lefx\Big].
\endaligned
\tag{9}
$$

The left-hand side depends only on $M^2$ so we can define the deformed mass 
square Casimir as
$$
\aligned
\mti^2 & = g^{00}  \Big(2\kappa \sh \kapti\Big)^2  + 4\kappa g^{0l} 
\pti_l \lefx \sh \kaptid\\
& + g^{rs} \pti_r \lefx \pti_s \lefx.
\endaligned
\tag{10}
$$
Obviously, $\mti^2 \to M^2 $ as $\kappa \to \ft$. he following relation 
between $M^2$ and  $\mti^2$ arises as a consequence of (6) and (9)
$$
A^2\Big( \dfrac{4\kappa^2}{\mti^2} + g_{00}\Big) = M^2.
\tag{11}
$$
In principle, we can also find the second Casimir operator --- the 
Pauli-Lubanski invariant. However, the resulting formula is quite involved 
and will be not quoted here.

Now let us consider the case of Weyl \4. As we stressed before, it is 
consistent only provided $g_{00} = 0$. Since the \8 is a sub\4 of \9, our 
formulae (7) are still valid. However, $A$ and $C$ should be subjected to 
some further conditions providing the proper action of dilatation operator 
$D$ on new momenta. It  appears, as a result of simple computation, that $C$ 
should be $M^2$-independent constant. If we put $C = \kappa$, we obtain 
finally the following \0 map
$$
\aligned
& \pti_0 = \kappa \ln \Big( \dfrac{P_0 + \kappa}{\kappa} \Big),\\
& \pti_i  = \dfrac{\kappa}{P_0 + \kappa} P_i  + \dfrac{M^2}{2(P_0 + 
\kappa)}  g_{i0}.
\endaligned
\tag{13}
$$

We have shown that formally the classical Poincar\'e (Weyl) \4 and the 
$\kappa$-Poincar\'e  ($\kappa$-Weyl) \4 are equivalent in the \4ic sector. 
This is no longer true in the co\4 sector: the standard Poincar\'e \4 is 
co\7ative while the $\kappa$-Poincar\'e one is not.

{\bf Acknowledgment.}
\flushpar
The numerous discussions with P. Kosi\'nski and P. Ma\'slanka are kindly 
acknowledged.

\enddocument